\documentclass{amsart}
\usepackage{a4wide}

\usepackage{amscd}
\usepackage{amssymb,amsmath}  
\usepackage{latexsym} 
\DeclareFontEncoding{OT2}{}{} 
\usepackage[english]{babel}

\newcommand{\textcyr}[1]{%
 {\fontencoding{OT2}\fontfamily{wncyr}\fontseries{m}\fontshape{n}\selectfont #1}}

\newcommand{\sha}{{\mbox{\textcyr{Sh}}}}
\newcommand{\q}{\mathfrak{q}}
\newcommand{\p}{\mathfrak{p}}

\newcommand{\Q}{\mathbf{Q}}

\newcommand{\F}{\mathbf{F}}
\newcommand{\C}{\mathbf{C}}

\newcommand{\Ps}{\mathbf{P}}

\newcommand{\ra}{\rightarrow}
\newcommand{\da}{\downarrow}

\newcommand{\str}{\mathcal{O}}

\newcommand{\M}{\mathcal{M}}

\newcommand{\ras}[1]{\stackrel{#1}{\rightarrow}}
\newcommand{\las}[1]{\stackrel{#1}{\leftarrow}}
\newcommand{\dra}{\swarrow}

\renewcommand{\bar}[1]{\overline{#1}}
\renewcommand{\phi}{\varphi}

    \newtheorem{Lem}{Lemma}[section]
    \newtheorem{Prop}[Lem]{Proposition}
    \newtheorem{Thm}[Lem]{Theorem}

   \theoremstyle{definition}
    \newtheorem{Def}[Lem]{Definition}
    \newtheorem{Con}[Lem]{Conjecture}
    \newtheorem{Cons}[Lem]{Construction}

    \newtheorem{Rem}[Lem]{Remark}

    \DeclareMathOperator{\Pic}{Pic}

\DeclareMathOperator{\rank}{rank}
\DeclareMathOperator{\tor}{tor}

\DeclareMathOperator{\cha}{char}
\DeclareMathOperator{\Br}{Br}

\DeclareMathOperator{\et}{\acute{e}t}

\author{Remke  Kloosterman}
\address{Department of Mathematics and Computer Science, University of Groningen, PO Box 800, 9700 AV  Groningen, The Netherlands}
\email{r.n.kloosterman@math.rug.nl}

\thanks{The author would like to thank Ronald van Luijk, Jasper Scholten and Jaap Top for many valuable discussions on this topic. The result of Section~\ref{cons} is taken from the author's PhD thesis \cite[Chapter 4]{Proefschrift}.}
\date{\today}
\title{Elliptic $K3$ surfaces with geometric Mordell-Weil rank 15}
\begin{document}
\begin{abstract}
We prove that the elliptic surface  $y^2=x^3+2(t^8+14t^4+1)x+4t^2(t^8+6t^4+1)$ has geometric Mordell-Weil rank 15. This completes a list of Kuwata, who gave explicit examples of elliptic $K3$-surfaces with geometric Mordell-Weil rank $0,1,\dots, 14, 16, 17,18$.
\end{abstract}
\maketitle

\section{Introduction}
The Mordell-Weil rank $r$ of a Jacobian elliptic surface $\pi:X\ra C$ is defined as the rank
of the group of sections of $\pi$. If $X$ is a $K3$ surface, then it
follows easily that $C=\Ps^1$. If one works over a field of
characteristic 0, then it is well known that $0\leq r\leq 18$. (In positive
characteristic we know that  $0\leq r \leq 20$.)

 By a result of Cox \cite{CoxK} there exists a Jacobian elliptic $K3$ surface
 defined over $\C$ with any given Mordell-Weil rank $r$, with  $r$ an integer, $0\leq
 r \leq 18$. Actually, using a similar reasoning as in \cite{CoxK} 
one can show there are infinitely
 many  $18-r$-dimensional families of Jacobian  elliptic $K3$ surfaces defined over $\C$, with Mordell-Weil rank $r$.
The examples constructed in the proof of Cox  are not explicit: the
 existence of such examples follows from properties of a so-called  period map.

Kuwata \cite{Kuw} has given a list of explicit Weierstrass equations for
 elliptic $K3$ surfaces defined over $\Q$ with Mordell-Weil rank $r$
 (over $\bar{\Q}$)  for any $r$ between 0 and 18, except for the case $r= 15$.

The aim of this paper is to complete this list by producing an explicit example of  an elliptic $K3$ surfaces with Mordell-Weil rank 15. This is achieved in two steps. In Section~\ref{cons} we prove the following:

\begin{Thm}\label{fifteenTHM} Let $K$ be an algebraically closed field, with
  $\cha(K)\neq 2,3$. Let $a,b,c\in K$. Let $E_{a,b,c}/K(s)$ be the  curve given
  by the Weierstrass equation
\[ y^2=x^3+A_{a,b,c}(s)x + B_{a,b,c}(s), \]
with \[A_{a,b,c}(s)=4a^3 b^3 ((b-a)cs^8+(2ac+2bc+4ab)s^4+(b-a)c\] and \[B_{a,b,c}=16 a^5b^5 s^2((b-a)s^8+2(b+a)s^4+ (b-a)).\]

For a general $(a,b,c)\in K^3$ this defines an elliptic $K3$ surface with 24 fibers of type $I_1$  and Mordell-Weil rank at least 15. In case $K=\C$ a generic member of this family has Mordell-Weil rank 15.
\end{Thm}

The strategy of our proof is the following. We start with a Jacobian  elliptic $K3$ surface $\pi: Y \ra \Ps^1$ such that there are exactly 15 components of fibers of $\pi$ not intersecting the zero-section. The Shioda-Tate formula (Theorem~\ref{ST}) implies that $\rho(Y):=\rank NS(Y)\geq 17$.

The particular examples $\pi:X\ra \Ps^1$ presented here, allow a degree 8 base-change of $\pi$ such that its associated relatively minimal model $\varphi: X \ra \Ps^1$ has only irreducible fibers and $X$ is a $K3$ surface as well. One can show in many ways (either in a direct and elaborate way or, if $K=\C$, using a powerful result in Hodge theory) that $\rho(X)=\rho(Y)$. The Shioda-Tate formula (Theorem~\ref{ST}) implies that the Mordell-Weil rank of $\varphi$ is at least 15.

In Section~\ref{exampl} we give an explicit example:

\begin{Thm}\label{examp} Assume that $K=\bar{\Q}$. The elliptic $K3$ surface $\pi: X \ra \Ps^1$ with Weierstrass equation
 \[y^2=x^3+2(t^8+14t^4+1)x+4t^2(t^8+6t^4+1) \]
has Mordell-Weil rank 15.\end{Thm}

The surface $X$ is isomorphic to the surface obtained by choosing $(a,b,c)=(2,4,2)$ in the equations of Theorem~\ref{fifteenTHM}.

 We use the following strategy to prove Theorem~\ref{examp}.
The proof of Theorem~\ref{fifteenTHM} yields that the Mordell-Weil rank of $\pi$ is at least 15. To prove equality we do the following:

 From the Shioda-Tate formula (Theorem~\ref{ST}) it follows that it suffices to prove $\rho(X)\leq 17$.  Since elliptic $K3$ surfaces over finite fields satisfy the Tate conjectures, one can determine $\rho(X \bmod q)$ from the characteristic polynomial $P_2(t)$ of Frobenius on $H^2_{\et}(X \bmod q,\Q_{\ell})$. The polynomial $P_2(t)$ can be easily determined using  the Lefschetz fixed point formula and counting the number of points on $X \bmod q$. To prove $\rho(X)\leq 17$ we find two prime numbers $p_1,p_2$ of good reduction such that the reduction $X \bmod p_1$ and $X \bmod p_2$ have Picard number 18. This is the best possible bound one can hope for by only considering $\rho(X\bmod q)$, since $\rho(X \bmod q)$ is even: from the fact $\rho(X)\leq \rho(X\bmod q)$ for any prime $q$ of good reduction it follows $\rho(X)\leq 18$.

We now use a refinement of the Tate conjectures to prove that $\rho(X)\leq 17$. Let $G_q$ be the Gram-matrix of the intersection pairing on $NS(X \bmod q)$. If $\rho(X)$ would be 18, then $\det(G_{p_1})$ and $\det(G_{p_2})$ would differ by a square. We use the Artin-Tate conjecture (which is both a refinement of and equivalent to the Tate conjecture, hence it holds for our $K3$ surfaces $X \bmod p_i$) to determine $\det(G_{p_1})$ and $\det(G_{p_2})$ up to squares. From this we deduce that $\rho(X) \leq 17$.  
\section{Definitions \and Notation}
\label{Prel}


\begin{Def}\label{defbas} An \emph{elliptic surface} is a triple
  $(\pi,X,C)$ with $X$ a smooth projective surface, $C$ a smooth projective curve, $\pi$ is a morphism
  $X\rightarrow C$, such that almost all fibers are irreducible genus
  1 curves and $X$ is relatively minimal, i.e., no fiber of $\pi$
  contains an irreducible rational curve $D$ with $D^2=-1$.

  We denote by $j(\pi): C \rightarrow \Ps^1$ the rational function such that $j(\pi)(P)$ equals the $j$-invariant of $\pi^{-1}(P)$, whenever $\pi^{-1}(P)$ is non-singular.

A \emph{Jacobian elliptic surface} is an elliptic surface together with a section $\sigma_0: C \rightarrow X$ to $\pi$.
  The set of sections of $\pi$ is an abelian group, with $\sigma_0$ as the identity element. Denote this group by $MW(\pi)$.


Let $NS(X)$ be the group of divisors on $X$ modulo algebraic equivalence, called the {\em N\'eron-Severi group} of $X$.
The {\em Picard number}  $\rho(X)$ is by definition the rank of the
  N\'eron-Severi group of $X$. 
\end{Def}








\begin{Rem}\label{singtabel} Suppose we are working over a field not of characteristic 2 or 3. If $P$ is a point on $C$, such that $\pi^{-1}(P)$ is
  singular then $j(\pi)(P)$ and $v_p(\Delta_p)$ behave as in Table~\ref{kod}.
For proofs of these facts see \cite[p. 150]{BPV}, \cite[Theorem IV.8.2]{Silv2}, \cite[p. 46]{Tate} or \cite[Lecture 1]{MiES}.
\end{Rem}

 \begin{table}[hbtp]\[\begin{array}{|c|c|c|c|}
\hline
\mbox{Kodaira type of fiber over }P & j(\pi)(P) & v_p(\Delta_p) & \mbox{ number of components}\\
\hline
I_0^*                  & \neq \infty & 6 &1\\
I_\nu \;(\nu>0) & \infty &\nu &\nu+1\\
I_\nu^* \;(\nu>0) & \infty&6+\nu  &\nu+5 \\
II     & 0  & 2& 1\\
IV &0& 4& 3\\
IV^* & 0 & 8& 7 \\
II^* & 0 & 10 &9\\
III &1728 & 3& 2\\
III^*              & 1728 & 9&8\\ \hline \end{array}\]\caption{Classification of singular fibers}\label{kod}\end{table}

\begin{Def} Let $X$ be a surface, let $C$ and $C_1$ be curves. Let $\varphi: X \rightarrow C$ and $f: C_1 \rightarrow C$ be two morphisms. Then we denote by $\widetilde{X \times_C C_1}$ the smooth, relatively minimal model of the  fiber product of $X$ and $C_1$ over $C$.
\end{Def}


Recall the following theorem.
\begin{Thm}[{Shioda-Tate (\cite[Theorem 1.3 \& Corollary
    5.3]{Sd})}]\label{ST}  Let $\pi:X\rightarrow C$ be a Jacobian
  elliptic surface, such that $\pi$ has at least one singular fiber. Then the N\'eron-Severi
  group of $X$ is generated by the classes of $\sigma_0(C)$, a
  non-singular fiber, the components of the singular fibers not
  intersecting $\sigma_0(C)$, and the generators of the Mordell-Weil
  group. Moreover, let $S$ be the set of points $P$ such that $\pi^{-1}(P)$ is singular. Let $m(P)$ be the number of irreducible components of $\pi^{-1}(P)$, then
  \[ \rho(X) =2+ \sum_{P \in S} (m(P)-1)+\rank(MW(\pi)) \]
\end{Thm}

The following result will be used several times. It is a direct
consequence of the Shioda-Tate formula.
\begin{Thm}[{\cite[Theorem 10.3]{Sd}}]\label{ratTHM} Let $\pi : X\ra
  \Ps^1$ be a rational Jacobian elliptic surface, then the rank of the Mordell-Weil group is 8 minus the number of irreducible components of singular fibers not intersecting the zero section. \end{Thm}

Given a Jacobian elliptic surface  $\pi: X \ra C$ over a field $K$, we can associate an
 elliptic curve in $\Ps^2_{K(C)}$  corresponding to the generic fiber
 of $\pi$. This induces a bijection between isomorphism classes of Jacobian
 elliptic surfaces and isomorphism classes of elliptic curves over $K(C)$.

Two elliptic curves $E_1$ and $E_2$ are isomorphic over $K(C)$ if and
only if $j(E_1)=j(E_2)$ and the quotients of the discriminants of $E_1/K(C)$ and $E_2/K(C)$ is a 12-th power (in $K(C)^*$).

Assume that $E_1$, $E_2$ are elliptic curves over $K(C)$ with
$j(E_1)=j(E_2)\neq 0,1728$. Then one shows easily that the quotients of the discriminants of $E_1$ and $E_2$ equals $u^6$ for some  $u\in K(C)^*$. Hence $E_1$ and
$E_2$ are isomorphic over $K(C)(\sqrt{u})$. We call $E_2$ the twist
of $E_1$ by $u$, denoted by $E_1^{(u)}$. Actually, we are not interested in the function $u$, but in the places at which the valuation of $u$ is odd.





\begin{Def} Let $\pi:X \ra C$ be a Jacobian elliptic surface. Fix $2n$ points $P_i \in C(\bar{K})$. Let $E/K(C)$ be the Weierstrass model of the generic fiber of $\pi$.

A Jacobian elliptic surface $\pi': X' \ra C$ is called a \emph{(quadratic) twist} of $\pi$ by $(P_1,\ldots,P_n)$ if the Weierstrass model of the generic fiber of $\pi'$ is isomorphic to $E^{(f)}$, where $E^{(f)}$ denotes the quadratic twist of $E$ by $f$ in the above mentioned sense
and $f\in K(C)$ is a function such that  $v_{P_i}(f) \equiv 1 \bmod 2$ and $v_Q(f) \equiv 0 \bmod 2$ for all $Q\not \in \{P_i\}$.
\end{Def}
If $K=\bar{K}$ then the existence of a twist of $\pi$ by $(P_1,\ldots,P_{2n})$ follows
directly from the fact that $\Pic^0(C)$ is 2-divisible.
Moreover, if  we fix $2n$ points $P_1,\ldots P_{2n}$ then there
exist precisely $2^{2g(C)}$ twists by $(P_i)_{i=1}^{2n}$.

If $P$ is one of the $2n$ distinguished points, then the fiber of $P$ changes in the following way (see \cite[V.4]{MiES}).
\[
I_\nu \leftrightarrow I^*_\nu \;(\nu \geq 0) \;\;\; \;\;
II \leftrightarrow    IV^* \;\;\;\;\;
III  \leftrightarrow  III^* \;\;\;\;\;
IV  \leftrightarrow  II^*
\]

Let $\pi: X\ra  C$ be a Jacobian elliptic surface, $P_1,\dots P_{2n}\in C$
points. Let $\tilde{\pi}: \tilde{X}\ra C$ be a twist by the $P_i$. Let
$\varphi: C_1\ra C$ be a double cover ramified at the $P_i$, such that the
minimal models of base-changing $\varphi$ and $\tilde{\varphi}$  by
$\pi$ are isomorphic. Denote this model by $\pi_1:X_1\ra C_1$.

Recall that 
\begin{eqnarray} \rank(MW(\pi_1)) = \rank(MW(\pi))+\rank(MW(\tilde{\pi})).\label{twteqn}\end{eqnarray} 
Moreover, the singular fibers change as follows
\[ \begin{array}{|l|cccc|}
\hline
\mbox{Fiber of } \pi  \mbox{ at } P_i                 & I_\nu \mbox{ or
  }I_\nu^*& II \mbox{ or } IV^*&
III  \mbox{ or } III^*& IV \mbox{ or } II^* \\

\hline
\mbox{Fiber of } \pi_1 \mbox{ at } \varphi^{-1}(P_i) & I_{2\nu}& IV &
I_0^* & IV^*\\\hline\end{array} \]
\section{Construction}\label{cons}
Let $K$ be an algebraically closed field of characteristic different
from 2 and 3.

Consider the following construction:

\begin{Cons}\label{constru} Let $\pi: X \rightarrow \Ps^1$ be a Jacobian
  elliptic surface whose singular fibers are three fibers of type $I_1$ and one fiber of type $III^*$. Let $f\in K(t)$ be a function of degree two, such that the fibers of $\pi$ over the critical values of $f$ are non-singular. 

Let $\alpha,\beta \in \Ps^1$ be the two distinct points such that
 $f(\alpha)=f(\beta)$ is the point which fiber is of type $III^*$. Let
 $g$ be a degree 4 cyclic covering, with only ramification over  $\alpha,
 \beta$. Let $\varphi:Y \rightarrow \Ps^1$ be the non-singular
 relatively minimal model of the fiber product $X \times_{\Ps^1}
 \Ps^1$ with respect to $\pi$ and $f \circ g : \Ps^1 \rightarrow
 \Ps^1$.
%
\end{Cons}

\begin{table}
\[\begin{array}{ccccccc}
        &&    &      &  Y   & \ras{\varphi} &  \Ps^1 \\
      &&      & \dra & \da &      & \da g_2'   \\
\Ps^1 & \las{\tilde{\pi}_2} &\tilde{X}_2 &      &  X_2 & \ras{\pi_2} & \Ps^1  \\
    &&        & \dra &\da &     & \da g_2   \\
\Ps^1 & \las{\tilde{\pi}_1} &\tilde{X}_1 &      & X_1 & \ras{\pi_1} & \Ps^1  \\
  &&          & \dra &\da &     & \da f   \\
\Ps^1 & \las{\tilde{\pi}} &\tilde{X}   &      & X   & \ras{\pi} & \Ps^1  \\
\end{array} \]
\caption{Overview of all maps used in this section ($g=g_2\circ g_2'$).}
\end{table}

\begin{Prop} \label{conPrp} The Mordell-Weil rank of $\varphi$ (of
 Construction~\ref{constru}) is at least 15, and is precisely 15 if and only if the rank of the twist of $\pi$ by the two critical values of $f$ is 0.
\end{Prop}

\begin{proof} 

The assumptions imply that $X$ is a rational surface and hence using
Theorem~\ref{ratTHM} we have that $\rank MW(\pi)=1$. Let $\pi_1 : X_1 \rightarrow \Ps^1$ be 
the fiber product $X \times_{\Ps^1} \Ps^1$ with respect to $f  : \Ps^1 \rightarrow \Ps^1$ and $\pi$. Let $\tilde{\pi}: \tilde{X} \rightarrow \Ps^1$ be the twist of $\pi$ by the two critical values of $f$. Then by (\ref{twteqn}) and  Theorem~\ref{ratTHM}
\[ \rank(MW(\pi_1))=\rank(MW(\pi))+\rank(MW(\tilde{\pi}))=1+\rank(MW(\tilde{\pi})). \]
Note that $\pi_1$ has two fibers of type $III^*$ and six fibers of type $I_1$. Let $P_1$ and $P_2$ be the points with a fiber of type $III^*$.

Let $g_2 : \Ps^1 \rightarrow \Ps^1$ be the degree two function,
 with critical values $P_1$ and $P_2$. 
Define $\pi_2: X_2 \rightarrow \Ps^1$ to be the non-singular relatively
minimal model of the fiber product $X \times_{\Ps^1} \Ps^1$ with
respect to $\pi_1$ and $g_2$.

Let $\tilde{\pi}_1$ be the twist of $\pi_1$ by $P_1$ and $P_2$. Then $\tilde{\pi}_1$ has two fibers of type $III$ and 6 fibers of type $I_1$, hence the corresponding surface is rational and by Theorem~\ref{ratTHM} $\tilde{\pi}_1$ has Mordell-Weil rank 6. From this it follows that $\rank(MW(\pi_2))=7+\rank(MW(\tilde{\pi}))$. Furthermore, $\pi_2$ has two fibers of type $I_0^*$, and 12 fibers of type $I_1$.

Let $\tilde{\pi}_2$ be the twist of $\pi_2$ by the two points with fiber of type $I_0^*$. Then $\tilde{\pi}_2$ has 12 fibers of type $I_1$ and the corresponding surface is rational with Mordell-Weil rank 8. So
\[ \rank(MW(\varphi))=\rank(MW(\pi_2))+\rank(MW(\tilde{\pi}_2)) = 15 + \rank(MW(\tilde{\pi})). \]\end{proof}

%
%

\begin{Rem} If we suppose that $\rank(MW(\tilde{\pi}))=0$, then it is relatively easy to find explicit generators for $MW(\varphi)$. In that case the pull-backs of the generators of $MW(\pi)$, $MW(\tilde{\pi}_1)$, $MW(\tilde{\pi}_2)$ generate a subgroup of $MW(\varphi)$ of index $2^m$, for some $m\geq 0$.  Since all these three surfaces are rational, we can take a specific Weierstrass model for these surfaces such that all Mordell-Weil groups are generated by polynomials of degree at most 2. (See \cite{OS}.) \end{Rem}

\begin{Rem} In the case $K=\C$ there exists another proof. Since $Y$
  and $\tilde{X}$ are both $K3$ surfaces, and there exists a finite
  map between them, the Picard numbers of both surfaces coincide (see
  \cite[Corollary 1.2]{Ino}). From an easy exercise using Kodaira's classification of singular fibers it follows that the configuration of singular fibers of $\varphi$ is the one mentioned in the Theorem.
By Kodaira's classification of singular fibers and the Shioda-Tate formula~\ref{ST}  we conclude
\[2+15+\rank(MW(\tilde{\pi}))=\rho(X)=\rho(Y)=2+\rank(\varphi).\]
\end{Rem}

%

Proposition~\ref{conPrp} enables us to prove the first theorem.

\begin{proof}[Proof of Theorem~\ref{fifteenTHM}] Let $c \in K^*$ such that $c^2\neq -1$. Then the rational elliptic surface $E'_c$ associated to the Weierstrass equation 
\[y^2=x^3+t^3(t-c) x + t^5\]
has a fiber of type $III^*$ and three fibers of type $I_1$. 
One easily shows that if $E'_c\cong E'_{c'}$ then $c'^2=c^2$. 
(If $E'_c\cong E'_{c'}$ then  there exist be an automorphism $h: \Ps^1
\rightarrow \Ps^1$ fixing 0 and $\infty$, and a constant $\lambda
\in K$, verifying $h(t)^3(h(t)-c)=\lambda^4 t^3 (t-c')$ and
$h(t)^5=\lambda^6 t^5$. This implies that $\lambda^2=1$ and
$c'=\lambda c$.)


Let $a\neq b$ and 
\[ f_{a,b}(s)=\frac{4ab s}{(a-b)s^2-2(a+b)s+a-b}. \]
The critical values of $f_{a,b}$ are $a$ and $b$, and $f^{-1}(0)=\{0,\infty\}$. Hence by Proposition~\ref{conPrp} the elliptic surface associated to the Weierstrass equation 
\[y^2=x^3+f_{a,b}(s^4)^3(f_{a,b}(s^4)-c) x + f_{a,b}(s^4)^5\]
satisfies the properties stated in the theorem. After a coordinate
change, which clears denominators, we obtain the equation of $E_{a,b,c}$.

This family contains a three-dimensional sub-family of non-isomorphic elliptic surfaces, because it is a finite base change of a three-dimensional family of non-isomorphic elliptic surfaces.

Assume now that $K=\C$. Let $\M$ be the moduli space of Jacobian elliptic $K3$ surfaces (cf. \cite{MiMS}). Let $U\subset \M$ be the set of elliptic
surfaces with non-constant $j$-invariant. Let $NL_{18}$ be the locus in $\M$ corresponding to elliptic $K3$ with Picard number at least 18. 

Suppose that a generic twist of $E'_c$
would have positive Mordell-Weil rank. Then the constructed family
$E_{a,b,c}$ would map to a 3-dimensional component $C$ of $NL_{18}$,
moreover the general member of the family $E_{a,b,c}$ has
non-constant $j$-invariant, hence $\dim C \cap U=3$. From
\cite[Theorem 1.1]{NL} it follows that $\dim NL_{18}\cap
U\leq 2$, a contradiction. From Proposition~\ref{conPrp} it follows that
the generic member of $E_{a,b,c}$ has Mordell-Weil rank precisely 15.
%
%
\end{proof}

%

\section{A method for bounding the Picard number}\label{method}
In the previous sections we showed the existence of a family of $K3$ surfaces such that the general member has Mordell-Weil rank 15. In this section we give an explicit example. In general it is hard to determine the Mordell-Weil rank of a non-rational elliptic surface. In the case of elliptic $K3$ surfaces one might be able to compute the Mordell-Weil rank using the Tate conjectures (which are proven for elliptic $K3$ surfaces over finite fields.)

In this section $K$ is supposed to be a number field. Recall the following facts.

Suppose $\pi : Y \ra \Ps^1$ is an elliptic surface defined over a number field $K$. Fix a model for $\pi$ over $\str_K$, the ring of integers of $K$. Let $\p$ be a prime of  $\str_K$. Assume that $Y$ has good reduction at $\p$.
Let $q:=\# \str_K/\p \str_K$. Let $\bar{Y}$ be the reduction of $Y$ modulo $\p$. Then the reduction map
\[r_\p: NS(Y_{\bar{\Q}}) \ra NS(\bar{Y}_{\bar{\F_q}})\]
is injective (for example see \cite[Proposition 6.2]{RonaldHeron}). It turns out that one can determine the rank of $NS(\bar{Y}_{\bar{\F_q}})$ for varieties $\bar{Y}$ for which the Tate conjectures hold.
In any case, the below mentioned method gives us  an {\em upper} bound for $\rank NS(\bar{Y}_{\bar{\F_q}})$).


\begin{Def} Suppose $(\Lambda, \langle \cdot , \cdot \rangle)$ is a lattice. 
Let $G$ be a Gram matrix of $\Lambda$ with respect to $\langle \cdot, \cdot \rangle$. By definition,  the {\em discriminant} of $(\Lambda, \langle \cdot , \cdot \rangle)$ is the determinant of $G$, which we denote by $\Delta(\Lambda, \langle \cdot , \cdot \rangle)$ or $\Delta(\Lambda)$, if no confusion arises.\end{Def}

It is well-known that the N\'eron-Severi group (modulo torsion) together with the intersection pairing forms a lattice. This implies the following proposition.

\begin{Prop}\label{SurProp} Let $Y/K$ be a smooth projective surface.
Suppose $\p$ is a prime of good reduction. Let $q=\#\str_K/\p\str_K$. Suppose that the reduction map
\[ r_\p:NS(Y_{\bar{\Q}}) \otimes \Q \ra NS(\bar{Y}_{\bar{\F_q}})\otimes \Q\]
is an isomorphism. Then the determinant of the Gram matrices of the intersection pairings on $NS(Y_{\bar{\Q}})$ and $NS({\bar{Y}_{\bar{\F_q}}})$ differ by a square.
\end{Prop}

\begin{proof} Since $\p$ is a prime of good reduction we have that $NS(Y_{\bar{\Q}})$ is a sublattice of  $NS({\bar{Y}_{\bar{\F_q}}})$ (see \cite[Proposition 6.2]{RonaldHeron}). Our assumptions imply that both lattices have the same rank. A standard result in lattice theory gives
\[ \Delta(NS(Y_{\bar{\Q}}))= [NS(\bar{Y}_{\bar{\F_q}}):NS(Y_{\bar{\Q}})]^2 \Delta(NS(\bar{Y}_{\bar{\F_q}}))\]
where $\Delta(\Lambda)$ denotes the determinant of the Gram matrix of the lattice $\Lambda$. This yields the proposition.
\end{proof}

One can construct examples of surfaces such that for every prime $\p$ of good reduction the image of $r_\p$ is {\em not} of finite index in $NS(\bar{Y}_{\bar{\F_q}})$ (see \cite{ell}, \cite{ShiodaPicone}, \cite{tera} for examples such that $\rho(Y)=1$; an easy example is the Kummer surface $Y$ of $E\times E$, for an elliptic curve $E/K$ without potential complex multiplication, then $\rho(Y)=19$, while all good reductions have by Proposition~\ref{evenprop} an even Picard number). 
Proposition~\ref{SurProp} turns out to be useful in showing that $r_\p$ is not surjective. 
In Section~\ref{exampl} we give such an example.

\begin{Con}[Tate Conjecture]\label{tatecon} Let $Y/\F_q$ be a smooth surface. Let $F_q$ be the automorphism of $H^2_{\et}(Y,\Q_{\ell})$ induced by the Frobenius automorphism of $\F_q$. Let $Q(t)$ be $\det(I-tF_q \mid H^2_{\et}(Y,\Q_{\ell}))$. 
Then $\rho(Y)$ equals  the number of reciprocal zeroes of $Q$ of the from $q\zeta$, with $\zeta$ a root of unity.\end{Con}

This conjecture is known to be true for several classes of varieties, see for example \cite{tatconov}.

The following proposition shows a small set-back, namely that over finite fields the number of possible Picard numbers seems smaller than over fields of characteristic zero.

\begin{Prop}\label{evenprop} Let $Y/\F_q$ be a smooth projective surface for which Conjecture~\ref{tatecon} hold (e.g. $K3$ surfaces). Then $\rho(Y)- \dim_{\Q_\ell} H^2_{\et}(Y,\Q_{\ell})$ is even.\end{Prop} 

\begin{proof} After replacing $\F_q$ by a finite extension, if necessary, we may assume that we have a set of $\F_q$-rational divisors, generating  $NS(Y_{\bar{\F_q}})$. This implies that the characteristic polynomial of Frobenius on $H^2_{\et}(Y,\Q_{\ell})$ is of the form $g(t)(t-1/q)^{\rho(Y)}$, where $g\in \Q[t]$ is a polynomial such that all its reciprocal zeroes have absolute value $q$ \cite[Theorem 1.6]{DeligneWeil}. Since the Tate conjecture~\ref{tatecon} holds for $Y$, it follows that $g(\pm q)\neq 0$, hence $g$ has no zeroes on the real. In particular, $g$ has even degree.

From Conjecture~\ref{tatecon} it follows that 
\[ \dim_{\Q_{\ell}} H^2_{\et}(Y,\Q_{\ell})-\rho(Y)=\deg g(t), \]
which yields the Proposition.\end{proof}

\begin{Rem} In this remark we try to indicate why there might exist many surfaces such that $r_\p$ is not surjective for any prime of good reduction. 
Using the so-called period map one can show that every integer $r$, such that $1\leq r\leq 20$ occurs as the Picard number of an algebraic $K3$ surface $Y$ over $\C$. Since $\dim H^2_{\et}(Y,\Q_{\ell})=22$ it follows from Proposition~\ref{evenprop} that $\rho(\bar{Y}_{\bar{\F_q}})$ is even, hence one might expect many examples of $K3$ surfaces defined over number fields such that for every prime of good prime the image of the reduction map $r_\p$ is not of finite index in $NS(Y_{\bar{\F_q}})$.  \end{Rem}

Suppose we can show that for two different primes $\p_1,\p_2$ of good reduction  the rank of the N\'eron-Severi lattices is the same, but the discriminants of the N\'eron-Severi lattices differs by a non-square. Then we conclude by Proposition~\ref{SurProp} that the rank of $NS(Y_{\bar{\Q}})$ is at least one lower than the rank of $NS(\bar{Y}_{\bar{\F_{q_1}}})$. This method was suggested to the author by Ronald van Luijk, see also \cite{RonaldPicone}.

The above remarks are only useful, if for a given surface $Y/\F_q$ one can efficiently compute $\rho(Y)$ and the determinant of the N\'eron-Severi lattice of $Y$. In general this is not the case, but for surfaces for which Conjecture~\ref{tatecon} holds, this can be done. Milne \cite{MilneAT} proved that if $Y$ is a surface for which Conjecture~\ref{tatecon} holds then also the following conjecture holds:

\begin{Con}[Artin-Tate Conjecture]\label{conAT} Let $Y/\F_q$ be a smooth surface. Let $P_2(t)$ be the characteristic polynomial of the $q$-Frobenius on $H^2_{\et}(Y,\Q_{\ell})$. Let $F_q$ be the Frobenius automorphism of $\F_q$. Let $Q(t):=\det(1-tF_q | H^2_{\et}(Y,\Q_{\ell}))$. Then
\[\lim_{s\ra 1} \frac{Q(q^{-s})}{(1-q^{1-s})^{\rho'(Y)}} = \frac{(-1)^{\rho'(Y)-1} \# \Br(Y) \Delta(NS(Y_{{\F_q}}))}{q^{\alpha(Y)}( \# NS(Y_{\F_q})_{\tor})^2}, \tag{*}\]
with $\alpha(Y)=\chi(Y,\mathcal{O}_Y)-1+\dim \Pic^0(Y)$ and $\Br(Y)$ is the Brauer group of $Y$. With $NS(Y_{{\F_q}})$ we indicate the subgroup of $NS(Y_{\bar{\F_q}})$ generated by $\F_q$-rational divisors and with $\rho'(Y)=\rank NS(Y_{\bar{\F_q}})$.
\end{Con}

If $Y/\F_{q}$ is an elliptic $K3$ surface, such that $\gcd(q,6)=1$, then both the Tate and the Artin-Tate conjecture are known to be true (see \cite{NygaardOgus}). We are actually interested in the discriminant of $NS(Y_{\bar{\F_q}})$. This forces us to apply Conjecture~\ref{conAT} over a field extension such that $\rho'(Y)=\rho(Y)$. (For an overview of cases for which the Tate and the Artin-Tate conjecture holds, see \cite{tatconov}.)

The Tate conjecture reduces our problem to finding the characteristic polynomial of Frobenius on $H^2_{\et}(Y,\Q_{\ell})$. This can be done by using the Lefschetz fixed point formula, i.e., one calculates the trace of the Frobenius automorphism (and several of its powers) on the cohomology by counting points on the surface. From knowing these traces, one deduces the characteristic polynomial. (This is explained in detail in \cite[Section 7]{RonaldHeron}).

We discuss some of the other quantities that have to be computed. Since we are only interested in knowing $\Delta(NS(Y_{\bar{\F_q}})$ up to squares, we might just disregard all quantities in (*) that are a square.


\begin{Prop}\label{KwadProp} Suppose $q$ is a prime power, with $\gcd(q,6)=1$. Let $\pi: Y\ra \Ps^1$ be an elliptic $K3$ surface, defined over $\F_q$. Assume that $q$ is a square and $\rho(Y)=\rho'(Y)$. Then
\[\Delta (NS(Y_{\bar{\F_\q}}) ) \equiv - \lim_{s\ra 1} \frac{Q(q^{-s})}{(1-q^{1-s})^{\rho(Y)}} \bmod \Q^{*2}\]
\end{Prop}

\begin{proof}
It is known that for an elliptic surface $\pi: Y \ra \Ps^1$ that the Brauer group $\Br(Y)$ is isomorphic to $\sha(Y/\Ps^1)$ (\cite[Chapter 5]{CosDol}), where $\sha(Y/\Ps^1)$ is the Tate-Shafarevich group of $Y/\Ps^1$. It is classically known that the number of elements of $\sha(Y/\Ps^1)$ is a square (\cite[Remark 6.11]{MiDu}, although this Remark is not completely correct, it is correct in the case of an elliptic curve over a function field), so we may disregard $\#\Br(Y)$. The parity of $\rho'(Y)$ follows from
\[ \rho(Y)=\rho'(Y)\equiv \dim H^2_{\et}(Y,\Q_{\ell}) \equiv 0 \bmod 2,\]
using Proposition~\ref{evenprop}. Since $\gcd(q,6)=1$ we know that Conjecture~\ref{conAT} holds. Combining these facts yields the proposition.
\end{proof}

\section{Proof of Theorem~\ref{examp}}
\label{exampl}
We apply the above mentioned strategy in the following example.

Consider the following elliptic $K3$ surface $X$ associated to
\[ y^2=x^3-(2t-1)^3 (4t-1)^2x+t(2t-1)^3(4t-1)^3.\]
This surface has two fibers of type $I_0^*$ (at $t=1/4$ and $t=1/2$) and a fiber of type $III^*$ at $t=\infty$. One can easily show that 17 and 19  are primes of good reduction. The components of the singular fibers, the zero-section and the class of a smooth fiber generate a rank 17 sublattice of $NS(X)$. One can easily show that these generators considered over $\F_{17}$ (resp. $\F_{19}$) are rational over  $\F_{17^6}$ (resp. $\F_{19^6}$), this is a straight-forward application of Tate's algorithm \cite{Tate}. A more precise application of Tate's algorithm yields an explicit degree 17 factor $F_p$ of  the characteristic polynomial $Q_p$ of the Frobenius of $\F_p$ acting on $H^2_{\et}(\bar{X},\Q_{\ell})$. In particular, one obtains that all roots $\alpha_p$ of $F_p$ satisfy $\alpha_p^6=p^{-6}$, if $p\in \{17,19\}$. Since $\dim H^2_{\et}(\bar{X},\Q_\ell)$ equals 22, we have to find a complimentary degree 5 factor $\tilde{G}_p$ of $Q_p$. From the fact that all reciprocal roots of $Q_p$ have absolute value $p$, at least one of the roots of $\tilde{G}_p$ is $\pm 1/p$. By counting points over $\F_{17}$ and $\F_{17^2}$ and using Poincar\'e duality (cf. \cite[Section 7]{RonaldHeron}), we obtain that  $G_{17}=(17x-1)\tilde{G}_{17}$, with $\tilde{G}_{17}$ given by
\[ 1+17x+136x^2+4913x^3+83521x^4\]
and similarly for $p=19$, we obtain that $G_{19}=(19x+1)\tilde{G}_{19}$, with  $\tilde{G}_{19}$ given by
\[ 1-9x-228x^2-3249x^3+130321x^4.\]
One easily shows that both polynomials have no reciprocal root of the form $p\zeta$, with $\zeta$ a root of unity. This implies that $\rho(\bar{X}_{\bar{\F}_p})=18$, for $p=17,19$.


Let $H_p$ the polynomial obtained by taking all roots of $G_p$ to the power six. Then $(p^6X-1)^{18}H_p$ is the characteristic polynomial of the Frobenius of $\F_{p^6}$ acting on $H^2_{\et}(\bar{X}_{\F_{p^6}},\Q_{\ell})$ for $p=17,19$. 

Proposition~\ref{KwadProp} implies that
%
\[ \Delta (NS(X_{\F_{17^6}}))\equiv 5\cdot19\cdot101516605992547\cdot11\cdot875005421 \bmod \Q^{*2} \]
and
\[ \Delta (NS(X_{\F_{19^6}}))\equiv
809308043\cdot95814202607062823339\cdot2297\cdot774901\cdot7\cdot13\cdot419\cdot16620229 \bmod \Q^{*2}.\]

The strategy explained in the previous section now implies that $\rho(X)\leq 17$. From Proposition~\ref{conPrp} we know that $\rho(X)\geq 17$,  proving that $\rank NS(X)=17$, and $\rank MW(\pi)=0$.
Applying Proposition~\ref{conPrp} again gives that
 \[y^2=x^3+2(t^8+14t^4+1)x+4t^2(t^8+6t^4+1) \]
has Mordell-Weil rank 15.

\end{document}